\documentclass[journal, final, twocolumn, 11pt]{IEEEtran}
\usepackage{cite}
\usepackage[pdftex]{graphicx}
\graphicspath{{../pdf/}{../jpeg/}}
\DeclareGraphicsExtensions{.pdf,.jpeg,.png}
\usepackage{amsmath, amsfonts,amssymb, amsthm, bm}
\usepackage{booktabs}
\usepackage[para,online,flushleft]{threeparttable}
\usepackage{float}
\usepackage{hyperref}
\usepackage[caption=false, font=footnotesize]{subfig}
\DeclareMathOperator{\Tr}{Tr}
\newtheorem{theorem}{Theorem}[section]
\newtheorem{corollary}[theorem]{Corollary}
\newtheorem{lemma}[theorem]{Lemma}
\newtheorem{proposition}[theorem]{Proposition}
\newtheorem{observation}[theorem]{Observation}
\newtheorem{definition}{Definition}[section]

\begin{document}
\bstctlcite{IEEEexample:BSTcontrol}
\title{On the 1-Wasserstein Distance between Location-Scale Distributions and the Effect of Differential Privacy}
\author{Saurab~Chhachhi,~\IEEEmembership{Student Member,~IEEE,
} and Fei~Teng,~\IEEEmembership{Senior Member,~IEEE,}%
\thanks{Manuscript received \today; revised XXX. This work was supported by ESRC through  LISS DTP (ES/P000703/1).
(Corresponding author Saurab Chhachhi)}
\thanks{Saurab Chhachhi and Dr. Fei Teng are with the Control and Power Group, Department of Electrical and Electronic Engineering, Imperial College London, London, UK (e-mail: saurab.chhachhi11@imperial.ac.uk; f.teng@imperial.ac.uk)}
\thanks{Digital Object Identifier: XXX}}

\maketitle
\begin{abstract}
    We provide an exact expressions for the 1-Wasserstein distance between independent location-scale distributions. The expressions are represented using location and scale parameters and special functions such as the standard Gaussian CDF or the Gamma function. Specifically, we find that the 1-Wasserstein distance between independent univariate location-scale distributions is equivalent to the mean of a folded distribution within the same family whose underlying location and scale are equal to the difference of the locations and scales of the original distributions. A new linear upper bound on the 1-Wasserstein distance is presented and the asymptotic bounds of the 1-Wasserstein distance are detailed in the Gaussian case. The effect of differential privacy using the Laplace and Gaussian mechanisms on the 1-Wasserstein distance is studied using the closed-form expressions and bounds.
\end{abstract}
\begin{IEEEkeywords}
Wasserstein distance, location-scale families, closed form expression, analytical expression, differential privacy, statistical distance.
\end{IEEEkeywords}

\section{Introduction}
\subsection{Background}
The Wasserstein distance has been widely used as a metric to represent the distance between two probability measures. Uses range from loss functions for Generative Adversarial Networks (GAN)\cite{Arjovsky2017}, bounding generalisation errors of machine learning models\cite{Lopez2018}, determining estimator properties\cite{Wu2012} and even assessing data quality\cite{Ding2021}. One of the main advantages of the Wasserstein distance over other distances/divergences such as the Kullback-Liebler divergence (KLD) is that it is a metric. As such, it obeys four axioms: (1) identity of indiscernibles $d(x,x)=0$, (2) symmetry $d(x,y) = d(y,x)$, (3) triangle inequality $d(x,z) \leq d(x,y) + d(y,z)$, and  (4) non-negativity $d(x,y) \geq 0$ \cite{Angelis2021}. In addition, unlike the KLD, the Wasserstein distance is finite even when neither measure is absolutely continuous with respect to the other\cite{Cai2022}.

The analytical definition of the Wasserstein distance (also called the Kantorovich distance, Mallow's distance, $L_p$-metric and the Earth Mover's Distance for special cases) is given by\cite[Definition 6.1]{Villani2009}:
\begin{definition}
    ($p$-Wasserstein distance). Let $(X , d)$ be a Polish metric space, and let  $p \in [1,\infty)$. For any two marginal measures $\mu$ and $\nu$ on $X$ , the Wasserstein distance of order $p$ between $\mu$ and $\nu$ is given by:
    \begin{align}
        W_p(\mu,\nu) &= \left(\inf_{\pi \in\Pi(\mu,\nu)} \int_{\mathcal{X}} d(x,y)^p d\pi(x,y)\right)^{1/p}\\
        \begin{split}
            &= \inf\{\mathbb{E}[d(X,Y)^p]^{1/p}, \\ & \qquad \qquad \mu = F_{\mu}(X), \nu = F_{\nu}(Y) \} 
        \end{split}
        \label{eq:def_norm}
    \end{align}
where, $\Pi(\mu,\nu)$ denotes the collection of all measures on $\mathcal{X}$ with marginals $\mu$ and $\nu$. The set $\Pi(\mu,\nu)$ is also called the set of all couplings of $\mu$ and $\nu$.
\end{definition}

In general, the Wasserstein distance does not admit closed-form expressions. There are two exceptions: (1) when $X$ and $Y$ are Gaussian then $W^2_2(X,Y)$ admits a closed-form expression and (2) when the distributions are univariate $(d=1)$\cite{Panaretos2019}.

The 2-Wasserstein distance between two Gaussians, $X_1 \sim N(\mu_1, \Sigma_1)$ and $X_2 \sim N(\mu_2, \Sigma_2)$, has the following closed-form representation for distributions\cite[Theorem 2.2]{Takatsu2011}:
\begin{align}
    \begin{split}
        W_2^2(X_1,X_2) = \lVert \mu_1 &- \mu_2\rVert^2 + \Tr{\left(\Sigma_{1}\right)} + \Tr{\left(\Sigma_{2}\right)} \\
        & - 2\Tr{\left(\left(\Sigma_{1}^{1/2}\Sigma_{2}\Sigma_{1}^{1/2}\right)^{1/2}\right)}
    \end{split}
    \label{eq:w2}
\end{align}
where $\Tr(\Sigma_i)$ is the trace of the covariance matrix $\Sigma_i$.

Although the above analytical expression is valid in the special case of  Gaussians it  provides a lower bound on the 2-Wasserstein distance for any symmetric distribution with the same distributional parameters ($\mu_i, \Sigma_i)$ \cite[Theorem 2.5]{Albertos1996}. In the univariate case, this expression for the 2-Wasserstein distance is exact for any elliptical symmetric distribution\cite{Gelbrich1990}. Given the closed-form expression for the 2-Wasserstein distance it has been popular. However, depending on the application, the 1-Wasserstein distance may be more desirable. For example, assuming $d(x,y)$ is the Euclidean distance the 1-Wasserstein distance is more robust to outliers as compared to the 2-Wasserstein distance (similar to difference between the mean absolute loss and the mean squared loss in linear regression models). Additionally, the 1-Wasserstein distance is used extensively for GANs because of the properties of its dual representation (restriction to 1-Lipschitz functions) as a special case of the Kantorovich-Rubenstein duality theorem which does not extend to the 2-Wasserstein distance \cite{Arjovsky2017,Panaretos2019}. Although, we note that GANs based on the 2-Wasserstein distance have also been explored\cite{Zhu2020}.

In the univariate case, the Wasserstein distance simplifies to a function of the difference between the quantile functions ($F^{-1}_X(q) = \inf\{x: F_{X}(t) \geq q\}, q \in (0,1)$):
\begin{align}
    W_p(X,Y) = \left(\int_{0}^{1} \lvert F^{-1}_X(q) - F^{-1}_Y (q)\rvert^p dq \right)^{1/p}
    \label{eq:w_quant}
\end{align}

Additionally, when $p=1$, an alternative expression in terms of the cumulative distribution functions (CDF) can be obtained:
\begin{align}
    W_1(X,Y) = \int_{\mathbb{R}} \lvert F_X(t) - F_Y(t)\rvert dt 
    \label{eq:w1_cdf}
\end{align}
We note here that the above representation shows that in the univariate case the 1-Wasserstein distance is the area between the marginal CDFs. This observation also leads to another definition of the 1-Wasserstein distance based on copulas\cite{Angelis2021}.

\subsection{Motivation}
Although (\ref{eq:w1_cdf}) provides a practical method for calculating the univariate 1-Wasserstein distance in many applications, a closed-form or analytical expression directly in terms of distributional parameters remains desirable. Indeed a closed-form/analytical representation would be more computationally efficient and convenient as it does not require the evaluation of an integral\cite{Angelis2021} and bypasses the need to conduct Monte-Carlo simulations. A closed-form expression also allows for exact solutions to, for example, DRO problems\cite{Farokhi2022}. It would also provide exact expressions for the effect of noise addition in differentially-private data analysis when the Wasserstein distance may be used as a metric of data utility\cite{Chhachhi2021}. In this case it is also desireable to be able to compute the Wasserstein distance privately to avoid potential privacy infringements. 

In the univariate discrete case, the Wasserstein distance is equivalent to the cardinality of set intersection between the two histograms\cite{Justicia2020}. This can be calculated privately using multiparty computation mechanisms known as Private Set Intersection-Cardinality (PSI-CA) which have at least $O(n)$ complexity. Closed-form expressions would allow for more efficient private calculation mechansism. However, to the best of our knowledge, a closed-form representation, in terms of distributional parameters, for the 1-Wasserstein distance for many widely used distributions (e.g. Gaussians) is not available either in the multivariate case or the univariate case. 

\subsection{Contribution}
This paper will focus on the 1-Wasserstein distance between independent univariate distributions belonging to a location-scale family.

\begin{definition}
    (Location-scale Distribution) For  $\alpha \in \mathbb{R}$ and  $\beta \in (0,\infty)$, let  $X=\alpha+\beta Z$. The two-parameter family of distributions associated with $X$ is called the location-scale family associated with the given (standard) distribution of $Z \sim (0,1)$ if its CDF is a function only of $\frac{x-\alpha}{\beta}$:
    \begin{align}
        F_{X}(x\mid \alpha, \beta) =  F\left( \frac{x-\alpha}{\beta} \right)
    \end{align}
with the standard CDF defined as:
    \begin{align}
        F_{Z}(x) = \Phi_{Z}(x)
    \end{align}
    Consequently, its quantile function can be expressed as:
    \begin{align}
        F^{-1}_{X}(q\mid \alpha, \beta) =  \alpha + \beta \Phi_{Z}^{-1}(q)
    \end{align}
    where $\Phi_{Z}^{-1}(q)$ is the quantile function for the standard distribution $Z$. 
\end{definition}

 For the avoidance of confusion we specify that, $\alpha$ is the location parameter and $\beta$ the scale parameter where as $\mu$ denotes the usual mean, $\sigma$ the standard deviation, and $\Sigma$ the covariance matrix in the multivariate case. When referring to a specific random variable within the location-scale family, with a slight abuse of notation $\Phi_{D}(x)$ and $\Phi_{D}^{-1}(q)$ denote the standard CDF and quantile functions for a distribution $Z \sim D(0,1)$, respectively (e.g. $N$ for the Gaussian or $Lap$ for the Laplace). Lastly, we use the term closed-form throughout to mean either truly closed-form expressions (i.e. those based solely on elementary functions) or analytical expressions (i.e. those which include functions such as the Gamma function, $\Gamma(k)$, or the standard Gaussian CDF, $\Phi_N(x)$, which can be efficiently computed from lookup tables).

Using the above representations for univariate distributions, we provide improved distribution specific closed-form upper bounds and exact expressions for the 1-Wasserstein distance between two independent univariate location-scale distributions. The expressions are based solely on the location and scale parameters ($\alpha$, $\beta$) of the two distributions in question and the standard quantile function $\Phi_{Z}^{-1}(x)$. 

In addition, we apply these expressions to show how the 1-Wasserstein distance can be used as a metric of quality for differentially-private data. We endogenously incorporate the effect of additive noise mechanisms (the Laplace and Gaussian Mechanisms) for differential privacy in the 1-Wasserstein distance.

\section{Closed-form Bounds}
There are a number of well established bounds on the 1-Wasserstein distances for distributions with finite first and second moments. These results apply to the univariate location-scale distributions studied in the paper. We note that these bounds are conventionally presented in terms of means and standard deviations ($\mu, \sigma$) rather than location and scale ($\alpha, \beta$). 
\subsection{Existing Bounds}
\label{sec:ex_bounds}
Below we present the tightest existing bounds for location-scale distributions before providing a new upper bound for univarite location-scale distributions and discussing the conditions under which this new bound is tighter than the extant literature.
\subsubsection{Upper Bound}
\begin{lemma}
Given two univariate independent distributions $X_1 \sim (\mu_1, \sigma_1)$ and $X_2 \sim (\mu_2, \sigma_2)$ within a location-scale family, the 1-Wasserstein distance between them is upper bounded by: 
\begin{align}
    W_1^{UB2}(X_1,X_2) = \sqrt{\left( \mu_1 - \mu_2 \right)^2 + \left(\sigma_{1} - \sigma_{2} \right)^2}
    \label{eq:ub_w2}
\end{align}
    \label{lm:ub_w2}
\end{lemma}

\begin{proof}
Given two spherical distributions $X_1 \sim (\mu_1, \Sigma_1)$ and $X_2 \sim (\mu_2, \Sigma_2)$, where the marginal distributions are orthogonal (i.e. $\Sigma_i = \sigma_i^2I_d$) the 2-Wasserstein distance admits a closed-form:
\begin{align}
    W_2^2(X_1,X_2) = \left( \mu_0 - \mu_1 \right)^2 + d\left(\sigma_{1} - \sigma_{0} \right)^2
    \label{eq:w2_sph}
\end{align}
This is also known as the Frechet distance \cite{Frechet1957}. Next, note that $W_p \leq W_q$ for $p\leq q$ by H\"{o}lder's inequality\cite[Remark 6.6] {Villani2009} meaning:
   \begin{align}
        W_1(X_1,X_2) &\leq \sqrt{W_2^2(X_1,X_2)} \\
          &= \sqrt{\left( \mu_1 - \mu_2 \right)^2 + d\left(\sigma_{1} - \sigma_{2} \right)^2}
    \end{align}  
In the univariate case $d=1$, which concludes the proof. 
\end{proof}

\subsubsection{Lower Bound}
\begin{lemma}
    Given two distributions $X_1$ and $X_2$ with $E[X_1] = \mu_1$ and $E[X_2] = \mu_2$ the 1-Wasserstein distance between them is lower bounded by \cite{Danica2015}:
    \begin{align}
        W_1^{LB}(X_1,X_2) = \lvert \mu_1 - \mu_2 \rvert
    \end{align}
    \label{lm:lb}
\end{lemma}
\begin{proof}
The definition of Wasserstein distance in (\ref{eq:def_norm}) is based on the expected norm between the difference of $X_1$ and $X_2$. As the norm is a convex function we can apply Jensen's inequality and the linearity of expectation to obtain:
    \begin{align}
        W_1(X_1,X_2) &= \inf\mathbb{E}[\lvert X_1-X_2\rvert]\\
        &\geq \lvert\mathbb{E}[X_1-X_2] \rvert \\
        &= \lvert\mathbb{E}[X_1]-\mathbb{E}[X_2] \rvert\\
        &= \lvert \mu_1 - \mu_2 \rvert
    \end{align}
\end{proof}
This result also agree with analysis in \cite[Proposition 3.2]{Angelis2021} where when one distribution dominates the other ($X_1 \succ X_2$), for example if $\sigma_1 = \sigma_2$ and $\mu_1 > \mu_2$) then the 1-Wasserstein distance is:
\begin{align}
    W_1(X_1,X_2) &= \mathbb{E}[X_1] - \mathbb{E}[X_2 ]\\ 
    &= W_1^{LB}(X_1,X_2)
\end{align}

\subsection{Location-scale Distributions}
We now present a new upper bound on the 1-Wasserstein distance specific to location-scale distributions.
\begin{theorem}
    Given two univariate independent random variables $X_1 = \alpha_1 + \beta_1Z$ and $X_2 = \alpha_2 + \beta_2Z$ the 1-Wasserstein distance between them is upper bounded by:
    \begin{align}
        W_1^{UB}(X_1,X_2) &= \lvert\alpha_1 - \alpha_2\rvert + \mathbb{E}[\lvert Z \rvert] \lvert\beta_1-\beta_2\rvert
        \label{eq:ub}
    \end{align}
    where $Z\sim (0,1)$ is a standard distribution within the location-scale family.
    \label{th:ub}
\end{theorem} 
\begin{proof}
    \begin{align}
        W_1(X_1,X_2) &= \int_0^1\lvert F^{-1}_1(q)- F^{-1}_2(q)\rvert dq \\
        \begin{split}
         &= \int_0^1\lvert \left(\alpha_1 + \beta_1\Phi_{Z}^{-1}(q)\right) \\
         &\quad \qquad - \left(\alpha_2 + \beta_2\Phi_{Z}^{-1}(q)\right)\rvert dq 
        \end{split}
    \end{align}
 by the triangle inequality $(\lvert x + y \rvert \leq \lvert x \rvert + \lvert y \rvert)$:
    \begin{align}
        \begin{split}
             W_1(X_1,X_2) \leq &\lvert \alpha_1 - \alpha_2 \rvert \\
             & + \lvert \beta_1-\beta_2 \rvert \int_0^1\lvert\Phi_{Z}^{-1}(q)\rvert dq
        \end{split}
    \end{align}
then note that $\int_0^1\lvert\Phi_Z^{-1}(q)\rvert dq = \mathbb{E}[\lvert Z \rvert]$.
\end{proof}

The new upper bound in Theorem \ref{th:ub} is linear in distributional parameters, as opposed to the existing bound in Lemma \ref{lm:ub_w2}. As $X_1$ and $X_2$ are from the same location-scale family we can reframe the expressions such that $\mu_i = \alpha_i$ and $\sigma_i = \beta_i$. We note that when $\alpha_1 = \alpha_2 / \mu_1 = \mu_2$ then (\ref{eq:ub_w2}) reduces to $\lvert \sigma_1 - \sigma_2\rvert$ and (\ref{eq:ub}) reduces to $\mathbb{E}[\lvert Z \rvert] \lvert\beta_1-\beta_2\rvert$. Therefore if $\mathbb{E}[\lvert Z \rvert] \leq 1$, (\ref{eq:ub}) will provide a tighter upper bound.

\section{Exact Analytical Expressions}
The closed-form bounds described above provide good linear as well as quadratic approximations however it remains desirable to obtain an exact expression for the 1-Wasserstein distance. To this end we provide an exact expression by revisiting the proof of Theorem \ref{th:ub} and modifying one step providing the main result of this paper.

\begin{theorem}
Given two univariate independent random variables $X_1 = \alpha_1 + \beta_1Z$ and $X_2 = \alpha_2 + \beta_2Z$ the 1-Wasserstein distance between them is:
    \begin{align}
        \begin{split}
            W_1(X_1,X_2)&= E[\lvert Y \rvert]
        \end{split}
        \label{eq:exact}
    \end{align}
    \label{th:exact}
where $Y = \left(\alpha_1 - \alpha_2\right) + \left(\beta_1-\beta_2\right) Z$.
\end{theorem}
\begin{proof}
    \begin{align}
        W_1(X_1,X_2)&= \int_0^1 \lvert \alpha_1 - \alpha_2 + \left(\beta_1-\beta_2\right)\Phi_{Z}^{-1}(q)\rvert dq \\
        &= \int_0^1 F^{-1}_{\lvert Y \rvert} (q) dq
    \end{align}
Note that the integrand defines the quantile function, $F^{-1}_{\lvert Y \rvert}(q)$, of a folded/absolute value random variable distributed as $Z$. Specifically, the underlying random variable $Y = \left(\alpha_1 - \alpha_2\right)+ \left(\beta_1-\beta_2\right)Z$. Via the substitution $q = F_{\lvert Y \rvert} (x)$ and $dq = F^{'}_{\lvert Y \rvert}(x) dx$:
\begin{align}
    &= \int_{-\infty}^{\infty} xf_{\lvert Y \rvert} (x) dx \\
    &= E[\lvert Y \rvert]
\end{align}
\end{proof}
We note that Theorem \ref{th:exact} can be extended to $p$-Wasserstein distances, using (\ref{eq:w_quant}), with the resulting value being $E[\lvert Y \rvert^p]^{1/p}$. However, we choose to focus on the 1-Wasserstein distance as in most cases for $p>1$ this quantity is not known or difficult to compute .

As Theorem \ref{th:exact} provides an exact expression for the 1-Wasserstein distance it may seem superfluous to include the new upper bound presented in Theorem \ref{th:ub}. However, there are a number of cases where Theorem \ref{th:ub} is either more desirable or indeed the only usable expression. A closed-form expression from Theorem \ref{th:exact} is predicated on the existence of a closed-form for the mean of the folded/absolute value random variable. In some cases this is not available at all or not in the general case ($\alpha_y \in \mathbb{R}, \beta_y \in \mathbb{R}^{+}$). For example, the mean of the folded Student's t distribution only has a convenient closed-form when $\alpha_y = 0$ (see Appendix \ref{sec:app_dist_table}), as a result only Theorem \ref{th:ub} can be used instead when $\alpha_y \neq 0$. 

Additionally, as will be discussed in greater detail in the next section, the upper bound provided in Theorem \ref{th:ub} can provide a more useful functional form (e.g. linear for Gaussians) which is particularly useful in application such as optimisation.    

\begin{figure}
    \centering
    \subfloat[$X_1 = 5 + 5Z, X_2 = \alpha_2 + 3Z$]{
    \includegraphics[width=0.45\columnwidth]{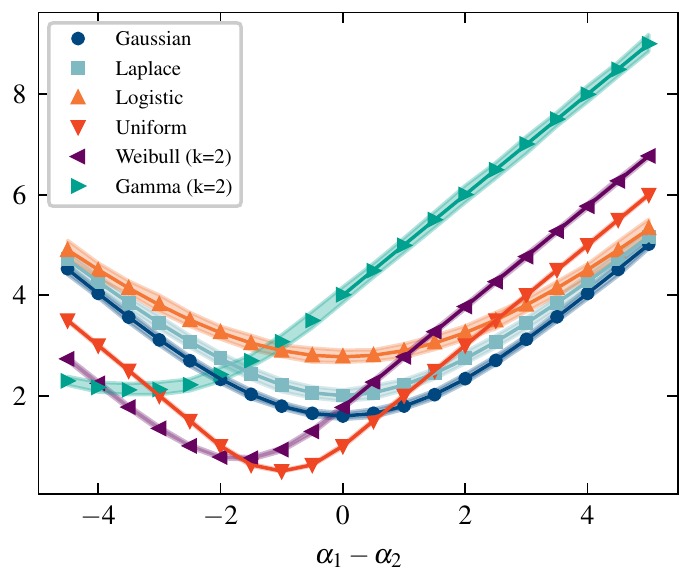}}
    \hfill
    \subfloat[$X_1=1 + 5Z, X_2 = \beta_2Z$]{
    \includegraphics[width=0.45\columnwidth]{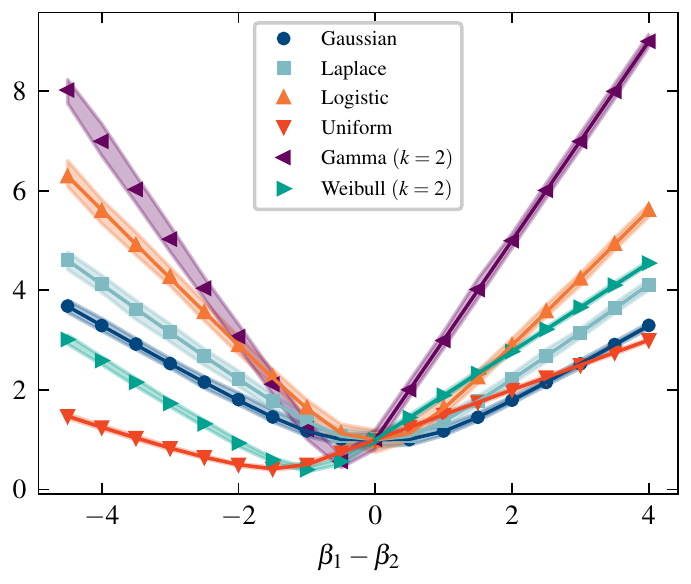}}
    \caption{1-Wasserstein distance for selected location-scale distributions}
    \label{fig:dist_bound}
\end{figure}

Figure \ref{fig:dist_bound} shows the 1-Wasserstein distance for selected distributions. The empirical distance was calculated using the \texttt{Python Optimal Transport (POT)} package\cite{flamary2021pot}. The empirical 1-Wasserstein distance (marker) is averaged over $N_r = 10^2$ simulations with $N_s = 10^4$ samples in each simulation. The shaded area indicates the 95\% confidence interval. The closed-form expression based on Theorem \ref{th:exact} are represented by the solid lines. We provide a list of closed-form expressions for the 1-Wasserstein distance between selected location-scale distributions in Appendix \ref{sec:app_dist_table}. For conventionally non-negative distributions such as the Weibull distribution, the 1-Wasserstein distance is simply $E[|Y|] = E[Y]$. However, in the general case the mean of the folded variable is not readily available in the literature. As such, Figure \ref{fig:dist_bound} does not include the closed-form values (solid lines) when $\alpha_1 - \alpha_2 < 0$ or $\beta_1 - \beta_2 < 0$ for the Gamma or Weibull distributions. We show a basic example of how we can extend closed-form expressions for the uniform case in Appendix \ref{sec:app_dist_table}.

\section{Univariate Gaussians}
\label{sec:gauss}
In this section we study in greater detail the 1-Wasserstein distance between independent univariate Gaussians using the closed-form bounds and exact expressions developed in the previous sections. We note here that for Gaussians $\mu=\alpha$ and $\sigma = \beta$ so we parameterise them using $\mu$ and $\sigma$ throughout.
\subsection{Improved Upper Bound}
In Section \ref{sec:ex_bounds} we showed a generic upper bound based on the 2-Wasserstein distance. However, a specific upper bound for 1-Wasserstein distance between Gaussian distributions, based on distributional parameters, has also been developed\cite{Chafai2010}.

\begin{lemma}
    Given two independent multivariate Gaussians $X_1 \sim N(\mu_1, \Sigma_1)$ and $X_2 \sim N(\mu_2, \Sigma_2)$ an upper bound for the d-dimensional 1-Wasserstein distance is\cite[Lemma 2.4]{Chafai2010}:
\begin{align}
    \begin{split}
        W_1(X_1,X_2) \leq& \lvert \mu_1 - \mu_2  \rvert \\
            &+ \Bigg(
                \sum^{d}_{j=1}
            \Big\{
                \left(\sqrt{\lambda_{1,j}}- \sqrt{\lambda_{2,j}}\right)^2
            \\ &+ 2 \sqrt{\lambda_{1,j}\lambda_{2,j}}\left(1-v_{1,j}\cdot v_{2,j}
            \right)
            \Big\}
            \Bigg)^{1/2}
    \end{split}
    \label{eq:ub_ndim}
\end{align}
\label{lm:ub_w1}
\end{lemma}
where $\mu_i$ are the means, $\lambda_{i,j}$ is the ordered spectrum of the d-dimensional Gaussians and $v_{i,j}$ is the associated orthonormal basis of the eigenvectors.

In the univaiate case (\ref{eq:ub_ndim}) simplifies to (denoted $W_1^{UB1}(X_1,X_2)$):
\begin{align}
    W_1(X_1,X_2) & \leq \lvert\mu_1 - \mu_2 \rvert + \lvert\sigma_{1} - \sigma_{2} \rvert
    \label{eq:ub_old}
\end{align}
A proof is provided in Appendix \ref{sec:app_old_ub}. We can recover (\ref{eq:ub_old}) which by applying the Cauchy-Schwartz inequality to (\ref{eq:ub_w2}) $\left( \sum_i x_i^2 \leq \left(\sum_i x_i\right)^2 \right)$: 
\begin{align}
    W_1(X_1,X_2) & \leq \lvert\mu_1 - \mu_2 \rvert + \sqrt{d}\lvert\sigma_{1} - \sigma_{2} \rvert\\
      & = \lvert\mu_1 - \mu_2 \rvert + \lvert\sigma_{1} - \sigma_{2} \rvert
\end{align}
We show that our new bound in Theorem \ref{th:ub} is tighter bound than either $W^{UB1}_1(X_1,X_2)$ or $W^{UB2}_1(X_1,X_2)$ under certain conditions. Although $W^{UB2}_1(X_1,X_2)$ provides a tighter bound in general, it is possible to obtain a linear bound that is always tighter than $W^{UB1}_1(X_1,X_2)$ and also tighter than $W^{UB2}_1(X_1,X_2)$ when $\mu_1 = \mu_2$.
\begin{corollary}
Given two univariate independent Gaussians $X_1 = N(\mu_1,\sigma_1^2)$ and $X_2 = N(\mu_2,\sigma_2^2)$ the 1-Wasserstein distance is upper bounded by:
    \begin{align}
            W_1(X_1,X_2)&\leq \lvert\mu_y\rvert + \sqrt{\frac{2}{\pi}}\lvert\sigma_y\rvert 
        \label{eq:ub_gauss}
    \end{align}
where $\mu_y = \mu_1 - \mu_2$ and $\sigma_y = \sigma_1 - \sigma_2$
\end{corollary}

\subsection{Exact Analytical Expression}
Using Theorem \ref{th:exact}, the 1-Wasserstein distance between two univarite Gaussians is the mean of a folded Gaussian\cite[Equation 7]{Tsagris2014}.
\begin{corollary}
Given two univariate independent Gaussians $X_1 \sim N(\mu_1, \sigma_1^2)$ and $X_2 \sim N(\mu_2, \sigma_2^2)$ the 1-Wasserstein distance is equal to the mean of a folded Gaussian $E[|Y|]$ where $Y \sim N(\mu_y = \mu_1 - \mu_2, \sigma^2_y = (\sigma_1-\sigma_2)^2)$:
    \begin{align}
        \begin{split}
            W_1(X_1,X_2)= \lvert\mu_y&\rvert \left[ 1- 2 \Phi_N\left(-\frac{\lvert\mu_y\rvert}{\lvert\sigma_y\rvert}\right) \right] \\ 
            &  + \lvert\sigma_y\rvert \sqrt{\frac{2}{\pi}} \exp{\left(-\frac{\mu_y^2}{2\sigma_y^2}\right)}
        \end{split}
        \label{eq:ae_gauss}
    \end{align}
    \label{cor:gauss_exact}
\end{corollary}
As shown above the 1-Wasserstein distance can be expressed as a function of distributional parameters and the standard normal CDF, $\Phi_N(x)$.
\subsection{Asymptotic Bounds}
Given the exact analytical representation of the 1-Wasserstein distance in terms of the distributional parameters we determine the tightness of the closed-form upper and lower bounds and establish asymptotic bounds. By taking limits over the distributional parameters $(\mu_y =  \mu_1 - \mu_2, \sigma_y =  \sigma_1 - \sigma_2)$ we produce the following proposition.
\begin{proposition}
    Given two univariate independent Gaussians $X_1 \sim N(\mu_1, \sigma_1^2)$ and $X_2 \sim N(\mu_2, \sigma_2^2)$ the 1-Wasserstein distance between them converges asymptotically to:
    \begin{align}
        \lim_{\sigma_y \to 0 | \mu_y \to \infty/-\infty} W_1(X_1,X_2) &= W_1^{LB}(X_1,X_2) \\
        \lim_{ \sigma_y  \to \infty |  \mu_y  \to 0 } W_1(X_1,X_2) &= W_1^{Lim}(X_1,X_2)
    \end{align}
where $W_1^{Lim}(X_1,X_2) = W_1^{UB}(X_1,X_2) - \lvert \mu_1 - \mu_2 \rvert =  \sqrt{\frac{2}{\pi}} \lvert \sigma_1 - \sigma_2 \rvert$.
    \label{th:asymp}
\end{proposition}
\begin{proof}
    \begin{align}
        \begin{split}
        \lim_{ \lvert\sigma_y\rvert  \to 0} W_1(X_1,X_2) 
            &= \lvert\mu_y\rvert \left[ 1- 2 \Phi\left(-\frac{\lvert\mu_y\rvert}{0}\right) \right] \\ 
            &\qquad + (0) \sqrt{\frac{2}{\pi}} \exp{\left(-\frac{\mu_y^2}{2 (0)}\right)}
        \end{split}\\
        &= \lvert\mu_y \rvert \left[ 1- 2 (0) \right]\\
        &= \lvert \mu_1 - \mu_2 \rvert\\
        &= W_1^{LB}(X_1,X_2)
    \end{align}
\end{proof}
The proofs for the remaining limits ($\sigma_y \to \infty, \mu_y  \to 0, \mu_y  \to \infty/-\infty$) are provided in Appendix \ref{sec:app_asymp}. We note that the results in the degenerate cases have been reported in \cite[Example 2.5 \& 2.6]{Chafai2010}.

\subsection{Improved Lower Bound}
Based on the asymptotic analysis in Proposition \ref{th:asymp} we can see that a tighter lower bound can be obtained for univariate Gaussians. 

\begin{proposition}
    Given two univariate independent Gaussians $X_1 = N(\mu_1,\sigma_1^2)$ and $X_2 = N(\mu_2,\sigma_2^2)$ the 1-Wasserstein distance is lower bounded by:
    \begin{align}
            W_1(X_1,X_2)&\geq \max\left(\lvert \sigma_y \rvert \sqrt{\frac{2}{\pi}}, \lvert \mu_y \rvert  \right)
        \label{eq:lb_gauss}
    \end{align}
    \label{pr:lb_gauss}
\end{proposition}
\begin{proof}
Staring from (\ref{eq:ae_gauss}) we see that each component is positive. As a result:
    \begin{align}
        W_1(X_1,X_2)&\geq \lvert \sigma_y \rvert \sqrt{\frac{2}{\pi}} \exp{\left(-\frac{\mu_y^2}{2  \sigma_y^2}\right)} \\
        &\geq W_1^{Lim}(X_1,X_2)\\
        & = \lvert \sigma_y \rvert \sqrt{\frac{2}{\pi}}
    \end{align}
Additionally, we know from Lemma \ref{lm:lb} we have the lower bound $W_1^{LB}(X_1,X_2) = \lvert \mu_y \rvert $.    
\end{proof}
\begin{figure}
    \centering
    \subfloat[$X_1 \sim N(\mu_1,4), X_2 \sim N(5,9)$]{
    \includegraphics[width=0.45\columnwidth]{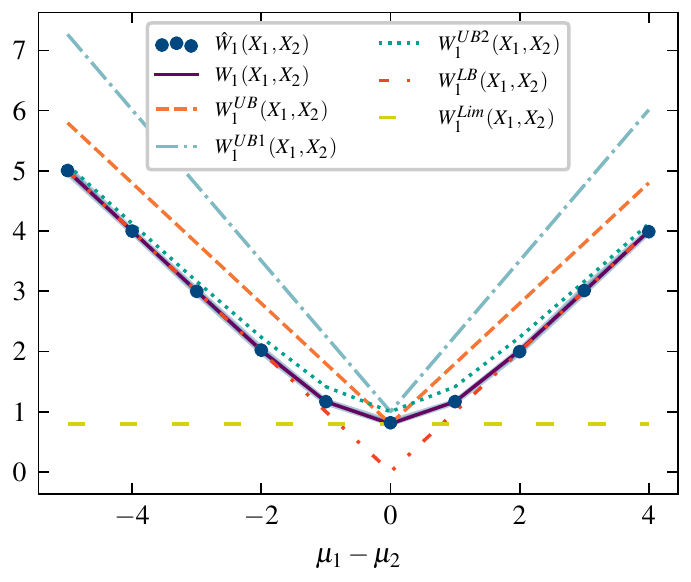}
    \label{fig:norm_a}}
    \hfill
    \subfloat[$X_1\sim N(2,\sigma_1^2), X_2\sim N(5,9)$]{
    \includegraphics[width=0.45\columnwidth]{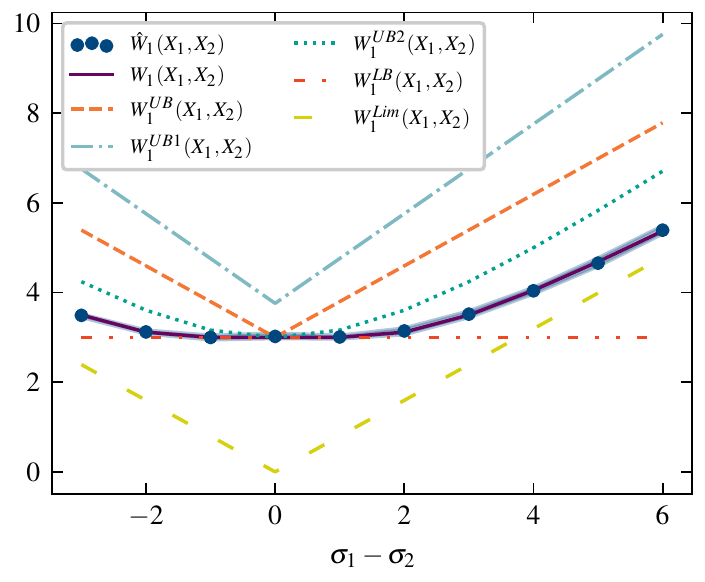}
    \label{fig:norm_b}}
    \caption{1-Wasserstein distance and bounds for univariate independent Gaussians}
    \label{fig:norm_bound}
\end{figure}

Figure \ref{fig:norm_bound} illustrates the improved bounds and exact expressions for the 1-Wasserstein distance between univariate Gaussians. $\hat{W_1}(X_1,X_2)$ is the empirical  1-Wasserstein distance averaged over $N_r = 10^2$ simulations with $N_s = 10^4$ samples in each simulation. The shaded area indicates the 95\% confidence interval. 

Although $W_1^{UB2}(X_1,X_2)$, defined in (\ref{eq:ub_w2}), is much tighter in general, the linear upper bound $W_1^{UB}(X_1,X_2)$ defined in (\ref{eq:ub_gauss}) is better than the existing linear upper bound $W_1^{UB1}(X_1,X_2)$ defined in (\ref{eq:ub_old}). We see that when either $\sigma_y \to 0$ or $\mu_y  \to \infty/-\infty$ the lower bound $W_1^{LB}(X_1,X_2)$ is tight. However, when either $\sigma_y \to \infty$ or $ \mu_y  \to 0$ they do not converge to the upper bound but rather an intermediate value $W_1^{Lim}(X_1,X_2)$, as discussed in Proposition \ref{pr:lb_gauss}. As a result the upper bound $W_1^{UB}(X_1,X_2)$ is only tight when $\sigma_y = 0$ resulting in $W_1^{UB}(X_1,X_2) = W_1^{LB}(X_1,X_2)$. This is clearly visible in Figure \ref{fig:norm_a}.

\section{Differential Privacy in the 1-Wasserstein Distance}
The Wasserstein distance is increasingly being used to characterise uncertainty, for example, in distributionally-robust optimisation (DRO)\cite{Farokhi2022} as well as a measure of data quality\cite{Ding2021}. Concurrently, there is a growing interest in ensuring user privacy for such data-driven applications. Differential privacy (DP), a privacy-preserving technique which achieves privacy through calibrated noise addition alters the data distribution and its resulting utility\cite{Chhachhi2021}. Incorporating this effect within the Wasserstein distance provides a unified analytical metric to assess intrinsic data quality and the utility degradation introduced by differentially-private noise.

 In this section we provide improved closed-form bounds for the 1-Wasserstein distance between a differentially-private data distribution $X_{DP} = X_1 + \text{DP-noise}$ and a reference (non-private) data distribution $X_2$. In particular, two popular noise addition mechanisms; (1) the Laplace mechanism and (2) the Gaussian mechanism, are studied. Additionally, in the case the data distributions $X_1$ and $X_2$ are Gaussian, a common assumption, we provide an exact expression and an approximation for the Gaussian and Laplace Mehcanisms respectively.

\subsection{Gaussian Mechanism}
The Gaussian Mechanism for differential privacy is defined as:
\begin{definition}
    (Gaussian Mechanism)  $\mathcal{M}_{N}(x,f,\epsilon, \delta)$ provides $(\epsilon,\delta)$-DP for a function $f(x)$\cite[Theorem A.1]{Dwork2014}:
    \begin{align}
        \mathcal{M}_{N}(x,f,\epsilon,\delta) = f(x) + X_N
    \end{align}
    where $X_N \sim N(0,\frac{2\ln(1.25/\delta)\Delta^2}{\epsilon^2}),  \Delta = \max\lvert f(x) = f(y) \rvert$, $\epsilon \in (0,1)$ is the privacy budget and $\delta$ is the probability of failure.
\end{definition}

Based on the definition above the differentially-private data distribution $X_{DP} = X_1 + X_N$. An upper bound
based on the triangle inequality and Jensen’s inequality
is provided in \cite[Theorem 2]{Farokhi2022}:
\begin{align} 
    \begin{split}
    W(X_1 +X_{N},X_2) &\leq W(X_1,X_2)\\ &\qquad +  \frac{\sqrt{2\ln(1.25/\delta)}\Delta}{\epsilon}
    \end{split}
    \label{eq:far2}
\end{align}

Using the corollaries outlined in Section \ref{sec:gauss} we generate improved upper bounds and in the case $X_1$ and $X_2$ are Gaussian, we can provide an exact expression for the 1-Wasserstein distance between differentially-private data and a reference distribution under the Gaussian mechanism. We first present the improved upper bound in Proposition \ref{pro:gauss}.

\begin{proposition}
    Given a $(\epsilon,\delta)$-DP data distribution $X_{DP} = X_1 + X_{G}$, where $X_N \sim N(0,\frac{2\ln(1.25/\delta)\Delta^2}{\epsilon^2})$, and a reference (non-private) data distribution $X_2$, the 1-Wasserstein between them is upper bounded by:
    \begin{align}
    \begin{split}
        W(X_1 + X_{N},X_2) &\leq W(X_1,X_2) \\
        &\qquad + \frac{2\Delta}{\epsilon}\sqrt{\frac{\ln(1.25/\delta)}{\pi}}
    \end{split}
    \end{align}
    \label{pro:gauss}
\end{proposition}
\begin{proof}
As the Wasserstein distance is a metric it obeys the triangle inequality:
    \begin{align} 
        W(X_1 +X_{N},X_2) &\leq W(X_1,X_2) + W(X_N,\mathfrak{\delta_0})
    \end{align}
where $\mathfrak{\delta_0}$ is the dirac delta distribution concentrated at 0. The second term can be reduced to:
\begin{align}
    W(X_N,\mathfrak{\delta_0}) &=  \sigma_N \sqrt{\frac{2}{\pi}} \quad \text{using (\ref{eq:exact})}  \\
    & = \frac{2\Delta}{\epsilon}\sqrt{\frac{\ln(1.25/\delta)}{\pi}}
\end{align}
\end{proof}

Next, if $X_1$ and $X_2$ are also Gaussian, we can provide an exact expression as $X_1$ and $X_N$ are independent and the the resulting differentially-private data is distributed as $X_{DP} \sim N(\mu_1,\sigma_1^2 + \frac{2\ln(1.25/\delta)\Delta^2}{\epsilon^2})$.
\begin{corollary}
Given a $(\epsilon,\delta)$-DP data distribution $X_{DP} = X_1 + X_{N}$, where $X_1 \sim N(\mu_1, \sigma_1^2$, $X_N \sim N(0,\frac{2\ln(1.25/\delta)\Delta^2}{\epsilon^2})$, and a reference (non-private) data distribution $X_2 \sim N(\mu_2, \sigma_2^2)$, the 1-Wasserstein between them is:
\begin{align}
    \begin{split}
        W(X_1 +X_{N},X_2) & =  \lvert\mu_y\rvert \left[ 1- 2 \Phi\left(-\frac{\lvert\mu_y\rvert}{\lvert\sigma_y\rvert}\right) \right] \\
        &\qquad + \lvert\sigma_y\rvert \sqrt{\frac{2}{\pi}} \exp{\left(-\frac{\mu_y^2}{2\sigma_y^2}\right)}
    \end{split}
\end{align}
where $\mu_y = \mu_1 - \mu2, \sigma_y = \sigma_{DP} - \sigma_2$ and $\sigma_{DP} = \sqrt{\sigma_1^2 + \frac{2\ln(1.25/\delta)\Delta^2}{\epsilon^2}}$.
\end{corollary}

\subsection{Laplace Mechanism}
The Laplace Mechanism for differential privacy is defined as:
\begin{definition}
    (Laplace Mechanism). $\mathcal{M}_{Lap}(x,f,\epsilon)$ provides $\epsilon$-DP for a function $f(x)$\cite[Definition 3.3]{Dwork2014}:
    \begin{align}
        \mathcal{M}_{Lap}(x,f,\epsilon) = f(x) + X_L
    \end{align}
    where $X_L \sim Lap(\alpha_l = 0, \beta_l = \frac{\Delta}{\epsilon}), \Delta = \max\lvert f(x) = f(y) \rvert$ and $\epsilon$ is the privacy budget.
\end{definition}
Based on the definition above the differentially-private data distribution $X_{DP} = X_1 + X_L$. A similar upper bound to (\ref{eq:far2}) for the Laplace mechanism is also provided in \cite[Theorem 2]{Farokhi2022}:
\begin{align}
    W(X_1 + X_{L},X_2) &\leq W(X_1,X_2) + \frac{\sqrt{2}\Delta}{\epsilon}
    \label{eq:far}
\end{align}
However, similar to Proposition \ref{pro:gauss}, in the univariate case it is possible to obtain a tighter upper bound.
\begin{proposition}
    Given a $\epsilon$-DP data distribution $X_{DP} = X_1 + X_{L}$, where $X_L \sim Lap(0, \frac{\Delta}{\epsilon})$, and a reference (non-private) data distribution $X_2$. The 1-Wasserstein between them is upper bounded by:
    \begin{align}
        W(X_1 + X_{L},X_2) \leq W(X_1,X_2) + \frac{\Delta}{\epsilon}
    \end{align}
    \label{prop:laplace1}
\end{proposition}
\begin{proof}
Again by the triangle inequality:
    \begin{align} 
        W(X_1 +X_{L},X_2) &\leq W(X_1,X_2) + W(X_L,\mathfrak{\delta_0})
    \end{align}
where $\mathfrak{\delta_0}$ is the dirac delta distribution concentrated at 0. The second term can be reduced to the following by applying Theorem \ref{th:exact} and the mean of the folded Laplace \cite[Proposition 2.3]{Liu2015}:
    \begin{align}
        W(X_L, \delta_0) &= \mathbb{E}[\lvert Y \rvert], \quad Y \sim Lap(0,\beta_y)\\
        &= \lvert 0 \rvert 
        + \lvert \beta_y \rvert \exp\left(-\frac{\lvert 0 \rvert}{\lvert \beta_y \rvert}\right) \\
        &= \lvert \beta_y \rvert \\
        &= \frac{\Delta}{\epsilon}
    \end{align}
\end{proof}

If $X_1$ and $X_2$ are Gaussian, we can provide an additional bound which is better than Proposition \ref{prop:laplace1} for larger $\epsilon$. The actual differentially-private data distribution will follow a Gaussian-Laplace distribution $X_{DP} \sim NL(\mu_1, \sigma_1, 1/b_l, 1/b_l)$ \cite{Reed2006}. The mean and variance are $E[X_{DP}] = \mu_1 $ and $Var[X_{DP}] = \sigma_1^2 + 2b_l^2$ respectively. Interestingly, for a given $b_l$ the Gaussian-Laplace is also a location-scale distribution. However, applying Theorem \ref{th:exact} would require the computation of the following quantity, $\int^{1}_{0} \lvert \beta_{DP}\Phi^{-1}_{NL}(q) - \sigma_{2}\Phi^{-1}_{N}(q)\rvert dq$ (where $\beta_{DP}(\sigma_1, b_l)$ is the scale parameter given that $X_{DP} = \mu_1 + \beta_{DP}Z_{NL}$), for which there is no closed form or convenient lookup table. Instead we observe that the Gaussian-Laplace can be approximated with high accuracy by a Gaussian with the same mean and variance. This is especially true when $b_l$ is smaller than or comparable to $\sigma_1$.
\begin{observation}
Given an $\epsilon$-DP data distribution $X_{DP} = X_1 + X_{L}$, where $X_1 \sim N(\mu_1, \sigma_1^2)$ and $X_L \sim Lap(0,\frac{\Delta}{\epsilon})$, and a reference (non-private) data distribution $X_2 \sim N(\mu_2, \sigma_2^2)$ the 1-Wasserstein between them can be approximated by the 1-Wasserstein between Gaussians. If $\sigma_1 \gtrapprox \sqrt{2}\frac{\Delta}{\epsilon}$ the following holds:
\begin{align}
    \begin{split}
         W(X_{DP}, X_2)  \approx \lvert\mu_y\rvert& \left[ 1- 2 \Phi_{N}\left(-\frac{\lvert\mu_y\rvert}{\lvert\sigma_y\rvert}\right) \right]\\ 
             &+ \lvert\sigma_y\rvert \sqrt{\frac{2}{\pi}} \exp{\left(-\frac{\mu_y^2}{2\sigma_y^2}\right)}
    \end{split}
    \label{eq:lap_app}
\end{align}
where $\mu_y = \mu_1 - \mu_2$, $\sigma_y = \tilde{\sigma}_{DP} - \sigma_2$, and $\tilde{\sigma}_{DP} = \sqrt{\sigma_1^2 + 2\frac{\Delta^2}{\epsilon^2}}$. 
\end{observation}

\begin{figure*}[htbp]
    \centering
    \subfloat[$X_1 \sim N(\mu_1,25), X_2 \sim N(5,4), \epsilon =  \frac{1}{4}$]{
    \includegraphics[width=0.3\textwidth]{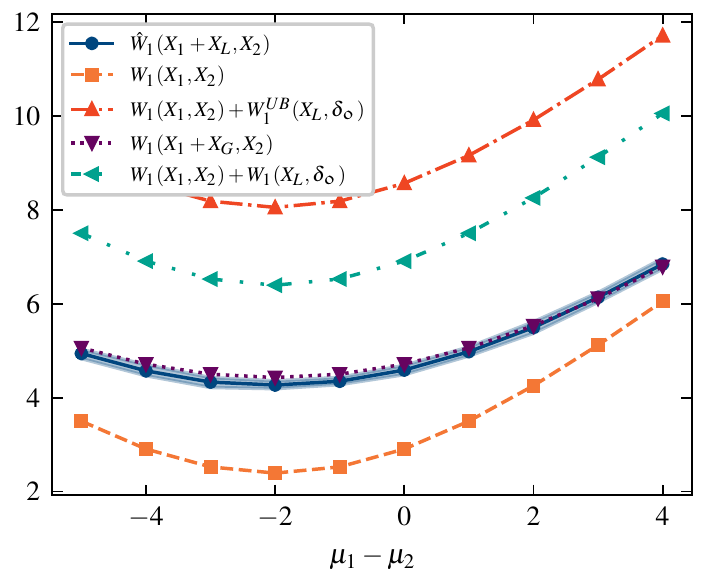}\label{fig:dp_a}}
    \hfill
    \subfloat[$X_1 \sim N(2,\sigma_1^2), X_2 \sim N(5,4),\epsilon =  \frac{1}{4}$]{
    \includegraphics[width=0.3\textwidth]{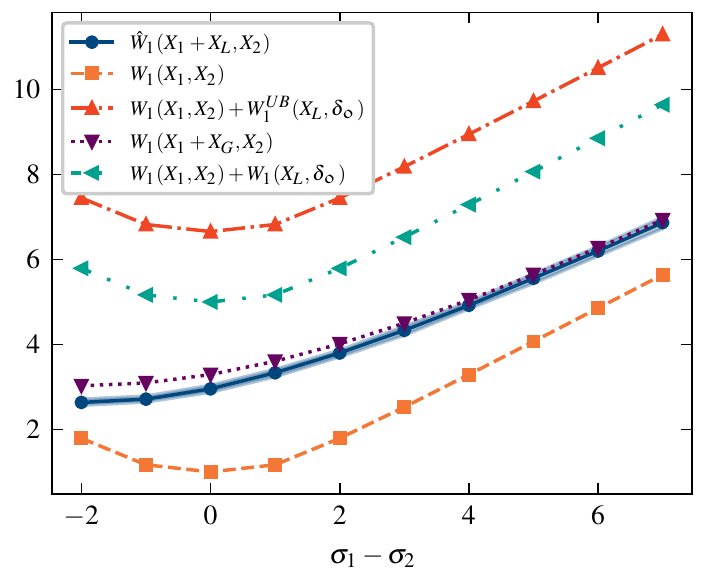}\label{fig:dp_b}}
    \hfill
    \subfloat[$X_1 \sim N(2,25), X_2 \sim N(5,4)$]{
    \includegraphics[width=0.3\textwidth]{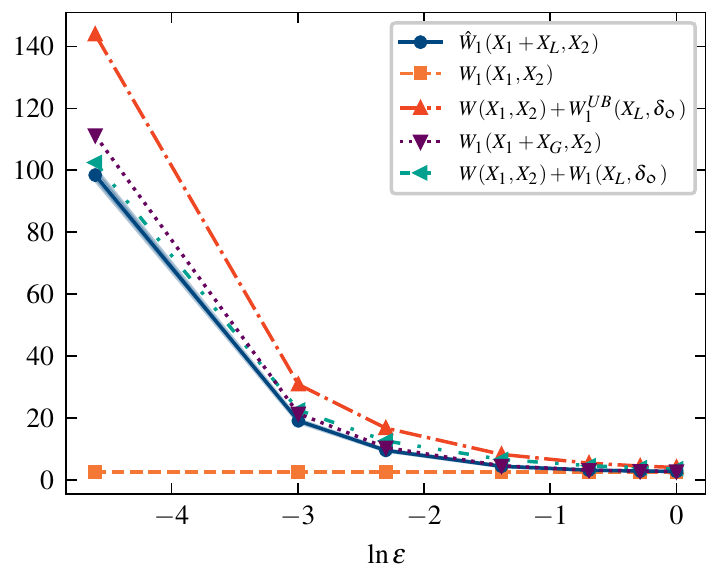}\label{fig:dp_c}}
    \\
    \subfloat[$X_1 \sim N(\mu_1,25), X_2 \sim N(5,4), \epsilon = 1$]{
    \includegraphics[width=0.3\textwidth]{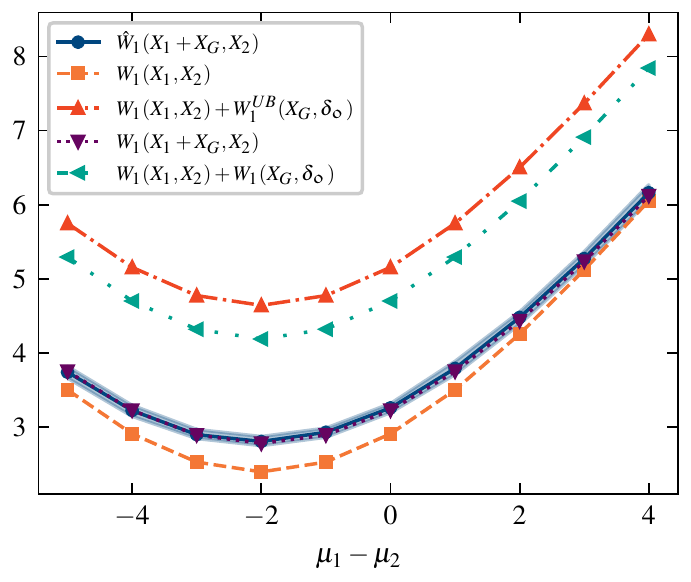}\label{fig:dp_d}}
    \hfill
    \subfloat[$X_1 \sim N(5,\sigma_1^2), X_2 \sim N(5,4),\epsilon = 1$]{
    \includegraphics[width=0.3\textwidth]{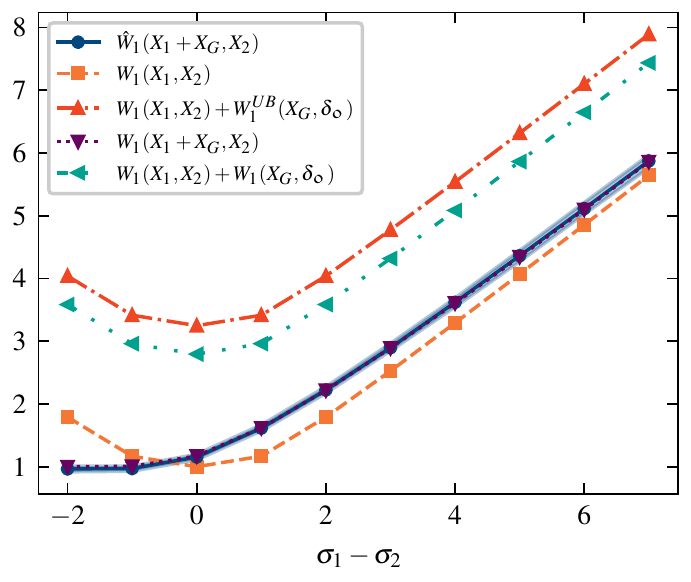}\label{fig:dp_e}}
    \hfill
    \subfloat[$X_1 \sim N(2,25), X_2 \sim N(5,4)$]{
    \includegraphics[width=0.3\textwidth]{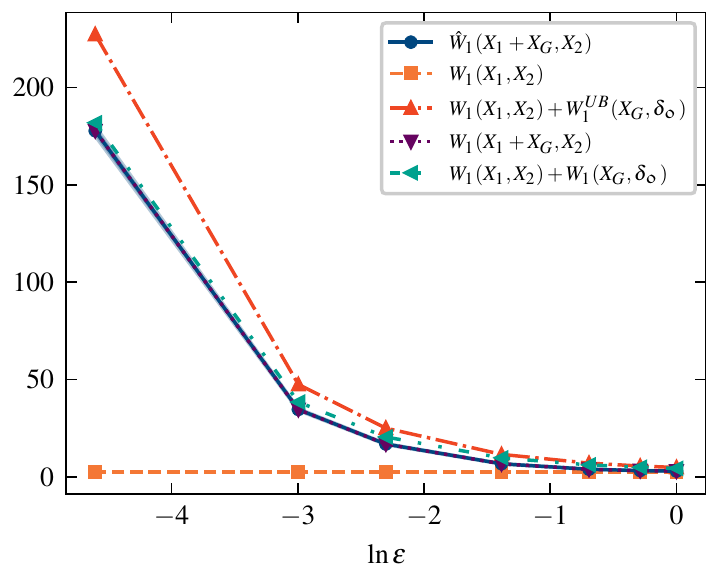}\label{fig:dp_f}}
    \caption{Differential Privacy in the 1-Wasserstein distance using Laplace and Gaussian Additive Noise Mechanisms}
    \label{fig:dp}
\end{figure*}

Figure \ref{fig:dp} illustrates the improved bounds provided above. We assume $\Delta = 1$ and $\delta = 10^-2$. The empirical 1-Wasserstein distance $(\hat{W}_1)$ averaged over $N_r = 10^2$ simulations with $N_s = 10^4$ samples in each simulation. The shaded area indicates the 95\% confidence interval. The bounds provided in (\ref{eq:far2}) and (\ref{eq:far}) are denoted $W_1(X_1,X_2) + W_1^{UB}(X_L, \mathfrak{\delta_0})$ and $W_1(X_1,X_2) + W_1^{UB}(X_G, \mathfrak{\delta_0})$ respectively. The improved bounds $W_1(X_1,X_2) + W_1(X_L, \mathfrak{\delta_0})$ and $W_1(X_1,X_2) + W_1(X_G, \mathfrak{\delta_0})$ perform significantly better than the previous bounds, especially for smaller privacy budgets($\epsilon$). The Gaussian approximation (\ref{eq:lap_app}) of the Laplace Mechanism denoted $W_1(X_1+X_G,X_2)$ in Figures \ref{fig:dp}\subref{fig:dp_a}-\subref{fig:dp_c} performs well compared to the upper bounds for larger privacy budgets($\epsilon$), where $\sigma_1 \gtrapprox \sqrt{2}\frac{\Delta}{\epsilon}$ .

\section{Conclusion}
This paper explored the properties of the 1-Wasserstein distance in the univariate case. We provided an exact analytical expression for the 1-Wasserstein distance between independent univariate location-scale distributions based solely on distributional parameters and special functions such as the standard Gaussian or Gamma function. In addition, a closed-form upper bound on the 1-Wasserstein distance for location-scale distributions is presented. In particular, for Gaussians this new bound is tighter than extant linear bounds and tighter overall when the means of the distributions are equal. The tightness of bounds was determined by exploring asymptotics of the exact analytical expression. Lastly, the effect of differentially-private noise addition on the 1-Wasserstein distance was investigated providing a tighter upper bound for both the Laplace and Gaussian mechanism and an exact expression in the case where the distributions are also Gaussian. Further work is needed to determine whether this approach can be extended to the multivariate case given that the theorems presented in this paper rely on the monotony of transport in the univariate case which allows the Wasserstein distance to be expressed in terms of quantile functions.

\bibliographystyle{IEEEtran}
\bibliography{wasserstein}

\appendices
\section{1-Wasserstein Distance between Selected Distributions}
\label{sec:app_dist_table}
Theorem \ref{th:exact} provides a convenient method for determining the 1-Wasserstein distance between two distributions within a location-scale family. For distributions with non-negative support (e.g. exponential, Weibull) the 1-Wasserstein is simply the mean of such a distribution ($Y$) with a location parameter of $\alpha_y = \left(\alpha_1 - \alpha_2 \right)$ and scale parameter of  $\beta_y = \left( \beta_1 - \beta_2 \right)$. For distributions with real support $(x \in \mathbb{R})$ (e.g. Gaussians) or bounded but both positive and negative support (e.g. $U(-1,1)$) the theorem requires an additional step to obtain a closed-form/analytical solution. The 1-Wasserstein distance is the mean of the absolute value (folded) of the distribution. In many cases explicit formulae for the folded distribution are readily available. Table \ref{tab:dist} summarises the 1-Wasserstein distance for widely used location-scale distributions. 

 Although uniformly distributed random variables are part of a location-scale family, a closed-form expression for the 1-Wasserstein distance between them does not have a single expression in the general case. Instead, there are two distinct cases as the mean of the resulting folded distribution is different depending on $\alpha_y$ and $\beta_y$.
 
 \begin{proposition}
     Given two univariate uniformly distributed random variables $X_1 = \alpha_1 + \beta_1Z$ and $X_2 = \alpha_2 + \beta_2Z$. The 1-Wasserstein distance between them is:
     \begin{align}
          W_1(X_1,X_2) &= \frac{1}{2}\left(\lvert a_y \rvert + \lvert b_y \rvert \right), \quad \text{ if } a_y,b_y \geq 0  \text{ or } \leq 0\\
             W_1(X_1,X_2) &= \frac{1}{2}\left(\frac{a_y^2  + b_y^2}{b_y  - a_y}\right), \quad \text{otherwise}
     \end{align}
     
 \end{proposition} 
 \begin{proof}
      By Theorem \ref{th:exact} the underlying random variable $Y \sim (\alpha_1 - \alpha_2) + (\beta_1 - \beta2)Z$. This is equivalent to a uniform random variable $Y \sim U(a_y = min((\alpha_1 - \alpha_2),(\beta_1 - \beta2)+(\alpha_1 - \alpha_2), b_y = max((\alpha_1 - \alpha_2),(\beta_1 - \beta2)+(\alpha_1 - \alpha_2))$. Below we derive the two cases depending on the resulting random variable $Y\sim U(a_y,b_y)$; (1) when $Y$ is non-negative or non-positive (i.e. $a_y,b_y \geq 0$ or $a_y,b_y \leq 0$) or (2) when $Y$ spans the origin (i.e. $a_y < 0, b_y > 0$).

\paragraph*{(1) $a_y,b_y \geq 0$ or $a_y,b_y \leq 0$}
In this case $\lvert Y \rvert \sim U(\lvert a_y \rvert, \lvert b_y \rvert)$ which means that $\mathbb{E}[\lvert Y \rvert] = \lvert\mathbb{E}[Y]\rvert$, resulting in a straight forward closed-form expression:
\begin{align}
    W_1(X_1,X_2) = \frac{1}{2}\left(\lvert a_y \rvert + \lvert b_y \rvert \right)
\end{align}

\paragraph*{(2) $a_y < 0, b_y > 0$}
In this case $\lvert Y \rvert$ will not have a uniform distribution, instead the negative support ($a_y \leq x \leq 0$) of $Y$ will be folded over into the positive domain. The expected value will then be:
\begin{align}
    E[\lvert Y \rvert] &= \int_{0}^{a_y} 2Cx dx + \int_{a_y}^{b_y} Cx dx\\
    &= \frac{C}{2} \left(a_y^2 + b_y^2 \right)
\end{align}
where $C$ is the normalising constant for $Y$:
\begin{align}
    2C a_y + C(b_y - a_y) = 1\\
    C = \frac{1}{b_y - a_y}
\end{align}
The 1-Wasserstein distance is then:

\begin{align}
    W_1(X_1,X_2) = \frac{1}{2}\left(\frac{a_y^2  + b_y^2}{b_y  - a_y} \right)
\end{align}
 \end{proof}

\begin{table}[htbp]
    \centering
     \caption{1-Wasserstein Distance for Selected Location-Scale Distributions}
    \label{tab:dist}
    \begin{threeparttable}
    \begin{tabular}{c c c}
    \toprule
        Distribution & $W_1(X_1,X_2)/\mathbb{E}[\lvert Y \rvert]$ \tnote{1} & Source\\
    \midrule
        Uniform\tnote{2} & $\displaystyle \begin{aligned}
        &\frac{1}{2}\left(\lvert a_y \rvert + \lvert b_y \rvert \right), a_y \text{ \& } b_y \geq 0 \text{ or } \leq 0 \\
        &\frac{1}{2}\left(\frac{a_y^2 + b_y^2}{b_y-a_y}\right), \quad \text{otherwise} 
        \end{aligned}$ &\\
         Gaussian &$\displaystyle \begin{aligned}
             \lvert\alpha_y\rvert& \left[ 1- 2 \Phi_{N}\left(-\frac{\lvert\alpha_y\rvert}{\lvert\beta_y\rvert}\right) \right]\\ 
             &+ \lvert\beta_y\rvert \sqrt{\frac{2}{\pi}} \exp{\left(-\frac{\alpha_y^2}{2\left(\beta_y\right)^2}\right)}
             \end{aligned}$
             & \cite{Tsagris2014}\\
         Laplace & $\displaystyle \begin{aligned}
         \lvert\alpha_y\rvert + \lvert\beta_y\rvert \exp{\left(-\frac{\lvert\alpha_y\rvert}{\lvert\beta_y\rvert}\right)}
         \end{aligned} $ &\cite{Liu2015}\\
         Logistic & $\displaystyle \begin{aligned}
         \lvert\alpha_y\rvert + 2\lvert\beta_y\rvert \ln{\left(1 + \exp{\left(-\frac{\lvert\alpha_y\rvert}{\lvert\beta_y\rvert}\right)} \right)}
         \end{aligned}$ &\cite{Cooray2006}\\
         Gamma \tnote{3}& $\displaystyle \alpha_y + k\beta_y, \quad \alpha_y,\beta_y \geq 0$ & \\
         Weibull \tnote{4}  & $\displaystyle  \alpha_y + \beta_y\Gamma(1+1/k), \quad \alpha_y,\beta_y \geq 0$   &\\
         Exponential\tnote{5} & $\displaystyle  \alpha_y + \beta_y, \quad \alpha_y,\beta_y \geq 0$     & \\
         Rayleigh\tnote{6} & $\displaystyle \alpha_y + \lvert\beta_y\rvert\sqrt{\frac{\pi}{2}}, \quad \alpha_y, \beta_y \geq 0$  &\\
         Student's t \tnote{7} & $\displaystyle \begin{aligned}
             2\lvert\beta_y\rvert\sqrt{\frac{\nu}{\pi}}\frac{\Gamma(\frac{\nu+1}{2})}{\Gamma(\frac{\nu}{2}) (\nu-1)}
         \end{aligned} $ & \cite{Psarakis1990} \\
    \bottomrule
    \end{tabular}
    \begin{tablenotes}
        \item[1] $Y \sim (\text{location: }\alpha_y = \left(\alpha_1 - \alpha_2\right), \text{scale: } \beta_y = \left(\beta_1 - \beta_2\right))$.\\
        \item[2] $Y \sim \alpha_y + \beta_y Z_{U} = U(a_y,a_y + b_y)$. The conventional upper and lower bounds of the uniform distribution $Y$ are $a_y = \min\left((\alpha_1 - \alpha_2), (\beta_1 - \beta_2) + (\alpha_1 - \alpha_2)\right) , b_y = \max\left((\alpha_1 - \alpha_2), (\beta_1 - \beta_2) + (\alpha_1 - \alpha_2)\right) $. \\
        \item[3] The Gamma distribution is a location-scale distribution for any given $k$. It is non-negative when $\alpha_y,\beta_y \geq0$ meaning $E[\lvert Y\rvert] = E[Y]$.\\
        \item[4] The Weibull distribution is a location-scale distribution for any given $k$. It is non-negative when $\alpha_y,\beta_y \geq0$ meaning $E[\lvert Y\rvert] = E[Y]$.\\
        \item[5] $\beta_y = \frac{ \lambda_1 - \lambda_2}{\lambda_1, \lambda_2}$ where $\lambda_i$ is the conventional inverse scale parameter. Equivalent to $Y \sim Weibull(\beta_y, k=1)$. Also reported independently in \cite{Chafai2010b} by evaluating (\ref{eq:w1_cdf}) directly.\\
        \item[6] Equivalent to $Y \sim Weibull(\sqrt{2}\beta_y, k=2)$.\\
        \item[7] Only applies for $\alpha_i=0$ and $\nu > 1$, where $\nu$ is the degrees of freedom.
    \end{tablenotes}
   \end{threeparttable}
\end{table}
\section{Simplification of Lemma \ref{lm:ub_w1}}
\label{sec:app_old_ub}
Below we provide a proof for the upper bound ($W_1^{UB1}(X_1,X_2)$) produced from Lemma \ref{lm:ub_w1}.
\begin{proof}
    The proof is provided by simplification from two dimensional case. Assuming the data is distributed normally with $X_1 \sim N(\mu_1, \Sigma_1)$ and $X_2 \sim N(\mu_2, \Sigma_2)$ their respective covariance matrix take the form:
\begin{align*}
    \Sigma_i = 
    \begin{bmatrix}
        {\sigma_{i,1}}^2 & \rho_i \\
        \rho_i & {\sigma_{i,2}}^2
    \end{bmatrix}
\end{align*}
The eigenvalues $\lambda_{i,j}$ of the covariance matrix can be determined by solving the following:
\begin{align*}
    \det\left(\Sigma_i - \lambda I_d\right)
     = \det
    \begin{bmatrix}
        {\sigma_{i,1}}^2 - \lambda_{i} & \rho_i \\
        \rho_i & {\sigma_{i,2}}^2 - \lambda_{i}
    \end{bmatrix}   = 0 
\end{align*}

\begin{align}
    \left({\sigma_{i,1}}^2 - \lambda_{i}\right)
    \left({\sigma_{i,2}}^2 - \lambda_{i}\right)
    - \rho_i^2 = 0 \\
    \lambda_n^2 - ({\sigma_{i,1}}^2 + {\sigma_{i,2}}^2)\lambda_i
    + \left( {\sigma_{i,1}}^2{\sigma_{i,2}}^2 - \rho_i^2 \right) = 0
\end{align}

The corresponding orthonormal eigenvectors $v_{i,j}$ can then be determined by solving the following:
\begin{align}
    \Sigma_i v_i = & \lambda_i \cdot v_i\\
    \begin{bmatrix}
        {\sigma_{i,1}}^2  & \rho_i \\
        \rho_i & {\sigma_{i,2}}^2
    \end{bmatrix} 
    \begin{bmatrix}
        v_{i,j,1} \\ v_{i,j,2}
    \end{bmatrix} 
    = &  \lambda_{i,j} \cdot
    \begin{bmatrix}
        v_{i,j,1} \\ v_{i,j,2}
    \end{bmatrix} 
\end{align}

If $\rho_i=0$, then:

\begin{align}
    \left({\sigma_{i,1}}^2 - \lambda_{i}\right)
    \left({\sigma_{i,2}}^2 - \lambda_{i}\right) = 0 \\
    \lambda_{i,1} = {\sigma_{i,1}}^2, \lambda_{i,2} = {\sigma_{i,2}}^2
\end{align}
The corresponding eigenvector for $\lambda_{i,1}$: 
\begin{align}
    \left( {\sigma_{i,1}}^2 - {\sigma_{i,1}}^2 \right) v_{i,1,1} = 0 \\
    \left( {\sigma_{i,2}}^2 - {\sigma_{i,1}}^2 \right) v_{i,1,2} = 0 \\
    v_{i,1} = \begin{bmatrix}
        1 \\ 0
    \end{bmatrix} 
\end{align}
and for $\lambda_{i,2}$: 
\begin{align}
    \left( {\sigma_{i,2}}^2 - {\sigma_{i,1}}^2 \right) v_{i,2,1} = 0 \\
    \left( {\sigma_{i,2}}^2 - {\sigma_{i,2}}^2 \right) v_{i,2,2} = 0 \\
    v_{i,2} = \begin{bmatrix}
        0 \\ 1
    \end{bmatrix} 
\end{align}
In the one-dimensional case therefore, $\lambda_i = \sigma_i^2$ and $v_i = 1$. The Wasserstein distance thus simplifies to:
\begin{align}
    \begin{split}
        W_1(X_1,X_2) &\leq \lvert \mu_1 - \mu_2  \rvert \\
            &+ \sqrt{
                \left(\sqrt{\sigma_1^2}- \sqrt{\sigma_2^2}\right)^2
            + 2 \sqrt{\sigma_1^2\sigma_2^2}\left(0
            \right)
            } 
    \end{split}\\
    &= \lvert \mu_1 - \mu_2 \rvert 
    +  \lvert\sigma_1 - \sigma_2 \rvert  
\end{align}
\end{proof}

\section{Proof of Asymptotic Convergence}
\label{sec:app_asymp}
The remaining proofs for Theorem \ref{th:asymp} are provided below.
\begin{proof}
\begin{align}
    \begin{split}
        \lim_{\lvert\sigma_y\rvert \to \infty} W_1(X_1,X_2) &= \lvert \mu_y \rvert \left[ 1- 2 \Phi\left(-\frac{\lvert\mu_y\rvert}{\infty}\right) \right] \\
        & \qquad +  \lvert\sigma_y \rvert  \sqrt{\frac{2}{\pi}} \exp{\left(-\frac{\left(\mu_y\right)^2}{2\left(\infty\right)^2}\right)}
    \end{split}\\
    &= \lvert \mu_y \rvert\left[ 1- 2\frac{1}{2} \right] + \lvert\sigma_y \rvert \sqrt{\frac{2}{\pi}}(1)\\
    &= \sqrt{\frac{2}{\pi}} \lvert \sigma_1 - \sigma_2 \rvert
\end{align}

\begin{align}
    \begin{split}
        \lim_{\lvert \mu_y \rvert \to 0} W_1(X_1,X_2) &= (0) \left[ 1- 2 \Phi\left(-\frac{0}{\lvert\sigma_y \rvert}\right) \right] \\
        & \qquad + \lvert\sigma_y \rvert\sqrt{\frac{2}{\pi}} \exp{\left(-\frac{\left(0\right)^2}{2\left(\sigma_y \right)^2}\right)}
    \end{split}\\
    &= \lvert\sigma_y \rvert \sqrt{\frac{2}{\pi}} (1)\\
    &= \sqrt{\frac{2}{\pi}}\lvert \sigma_1 - \sigma_2 \rvert
\end{align}

\begin{align}
     \begin{split}
         \lim_{\lvert \mu_y \rvert \to \infty} W_1(X_1,X_2) &=  \lvert \mu_y \rvert  \left[ 1- 2 \Phi\left(-\frac{\infty}{\sigma_y}\right) \right] \\
        & \qquad +  \sigma_y  \sqrt{\frac{2}{\pi}} \exp{\left(-\frac{\left(\infty\right)^2}{2\left(\sigma_y\right)^2}\right)}
    \end{split}\\
    &=  \lvert \mu_y \rvert  [1-(0)] +  \sigma_y \sqrt{\frac{2}{\pi}} (0)\\
    &=  \lvert \mu_y \rvert  \\
    &= W_1^{LB}(X_1,X_2)
\end{align}
\end{proof}
\end{document}